\documentclass{amsart}
\usepackage{amssymb,amsmath, amsthm,latexsym}
\usepackage{graphics}
\usepackage{amscd}
\usepackage{graphics}
\usepackage{here}
\newcommand{\cal}[1]{\mathcal{#1}}
\theoremstyle{plain}

\parskip=\bigskipamount

\let\egthree=\phi
\let\phi=\varphi
\let\varphi=\egthree

\begin{document}
\title{Bounded cohomology, cross ratios and cocycles}
\author{Ursula Hamenst\"adt}
\thanks
{Partially supported by Sonderforschungsbereich 611}
\date{August 24, 2005}

\begin{abstract}
We use cross ratios to describe second real continuous bounded
cohomology for locally compact $\sigma$-compact topological
groups. We also investigate the second continuous bounded
cohomology group of a
closed subgroup of the isometry group ${\rm
Iso}(X)$ of a proper hyperbolic geodesic metric space $X$ 
and
derive some rigidity results for ${\rm Iso}(X)$-valued cocycles.
\end{abstract}

\maketitle

\section{Introduction}

A geodesic metric space $X$ is called \emph{$\delta$-hyperbolic}
for some $\delta >0$ if it satisfies the $\delta$-thin triangle
condition: For every geodesic triangle in $X$ with sides $a,b,c$
the side $a$ is contained in the $\delta$-neighborhood of $b\cup
c$. If $X$ is proper then $X$ can naturally be compactified by
adding the \emph{Gromov boundary} $\partial X$.  The isometry
group ${\rm Iso}(X)$ of $X$, equipped with the compact open
topology, is a locally compact $\sigma$-compact topological group
which acts on $X\cup \partial X$ as a group of homeomorphisms. The
\emph{limit set} of a closed subgroup $G$ of ${\rm Iso}(X)$ is the
set of accumulation points in $\partial X$ of an orbit of the
action of $G$ on $X$; this limit set is a compact $G$-invariant
subset of $\partial X$. The group $G$ is called \emph{elementary}
if its limit set consists of at most two points.

Let $S$ be a standard Borel space and let $\mu$ be a Borel
probability measure on $S$. Let $\Gamma$ be a countable group
which admits a measure preserving ergodic action on $(S,\mu)$. An
${\rm Iso}(X)$-valued \emph{cocycle} for this action is a
measurable map $\alpha:\Gamma\times S\to {\rm Iso}(X)$ such that
$\alpha(gh,x)=\alpha(g,hx)\alpha(h,x)$ for all $g,h\in \Gamma$ and
$\mu$-almost every $x\in S$. The cocycle $\alpha$ is
\emph{cohomologous} to a cocycle $\beta:\Gamma\times S\to {\rm
Iso}(X)$ if there is a measurable map $\phi:S\to {\rm Iso}(X)$
such that $\phi(gx)\alpha(g,x)=\beta(g,x)\phi(x)$ for all $g\in
G$, $\mu$-almost every $x\in S$. Recall that a \emph{compact
extension} of a locally compact topological group $G$ is a locally
compact topological group $H$ which contains a normal compact
subgroup $K$ such that $G=H/K$ as topological groups. Extending
earlier results of Monod and Shalom \cite{MS04} we show in Section
4.

\bigskip

{\bf Theorem A:} {\it Let $G$ be a semi-simple Lie group with
finite center, no compact factors and of rank at least $2$. Let
$\Gamma<G$ be an irreducible lattice which admits a mixing measure
preserving action on a probability space $(S,\mu)$. Let $X$ be a
proper hyperbolic geodesic metric space and let
$\alpha:\Gamma\times S\to {\rm Iso}(X)$ be a cocycle; then one of
the following two possibilities holds.
\begin{enumerate}
\item $\alpha$ is cohomologous to a cocycle into
an elementary subgroup of ${\rm Iso}(X)$.
\item $\alpha$ is cohomologous to a cocycle
into a subgroup $H$ of ${\rm Iso}(X)$
which is a compact extension of a simple
Lie group $L$ of rank one,
and there is a continuous surjective homomorphism $G\to L$.
\end{enumerate}
}

\bigskip

The proof of this result uses second bounded cohomology for closed
subgroups of the isometry group ${\rm Iso}(X)$ of $X$. Here the
second bounded cohomology group of a locally compact topological
group $G$ with coeffients in a \emph{Banach module} $E$ for $G$,
i.e. a Banach space $E$ together with a continuous homomorphism of
$G$ into the group of linear isometries of $E$, is defined as
follows. For every $i\geq 1$, the group $G$ naturally acts on the
vector space $C_b(G^i,E)$ of continuous bounded maps $G^i\to E$.
If we denote by $C_b(G^i,E)^G\subset C_b(G^i,E)$ the linear
subspace of all $G$-invariant such maps, then the \emph{second
continuous bounded cohomology group} $H_{cb}^2(G,E)$ of $G$ with
coefficients $E$ is defined as the second cohomology group of the
complex
\begin{equation} 0\to
C_b(G,E)^G \xrightarrow{d} C_b
(G^2,E)^G\xrightarrow{d} \dots
\end{equation} with the usual
homogeneous coboundary operator $d$ (see \cite{M}). If $G$ is a
countable group with the discrete topology then we write
$H_b^2(G,E)=H_{cb}^2(G,E)$. Of particular importance is the case
that $E=\mathbb{R}$ with the trivial $G$-action which yields real
continuous second bounded cohomology. There is a natural map of
$H_{cb}^2(G,\mathbb{R})$ to the ordinary continuous second
cohomology group
$H_c^2(G,\mathbb{R})$ of $G$ which in general
is neither injective nor surjective.

If $G$ is a non-elementary closed subgroup of
the isometry group of a proper hyperbolic
geodesic metric space $X$ with limit set $\Lambda\subset
\partial X$ which does not act
transitively on the complement of the diagonal in $\Lambda\times
\Lambda$ then the kernel of the natural map
$H_{cb}^2(G,\mathbb{R})\to H_c^2(G,\mathbb{R})$ is infinite
dimensional \cite{H05b}. In contrast, we show in Section 3.

\bigskip

{\bf Theorem B:} {\it Let $X$ be a proper hyperbolic
geodesic metric space and let $G<{\rm Iso}(X)$ be
a non-elementary closed subgroup with limit set
$\Lambda$. If $G$ acts transitively on the complement of the
diagonal in $\Lambda\times \Lambda$ then the kernel
of the natural map $H_{cb}^2(G,\mathbb{R})\to
H_c^2(G,\mathbb{R})$ is trivial.}

\bigskip

In Section 2 we identify for an arbitrary locally compact
$\sigma$-compact topological group $G$ the group
$H_{cb}^2(G,\mathbb{R})$ with the vector space of measurable
anti-symmetric $G$-cross ratios on a strong boundary for $G$. As
an application we describe the kernel of the natural map
$H_b^2(\Gamma,\mathbb{R})\to H^2(\Gamma,\mathbb{R})$ for a lattice
$\Gamma$ in $PSL(2,\mathbb{R})$. Namely, recall that $\Gamma$ acts
as a group of orientation preserving isometries on the hyperbolic
plane ${\bf H}^2$ and hence it acts on the boundary
$S^1=\partial{\bf H}^2$ of ${\bf H}^2$ preserving the measure
class of the Lebesgue measure $\lambda$. We show that the kernel
of the natural map $H^2_b(\Gamma,\mathbb{R})\to
H^2(\Gamma,\mathbb{R})$ is naturally isomorphic to the vector
space of all finite $\Gamma$-invariant reflection-anti-invariant
finitely additive signed measures $\mu$ on the complement of the
diagonal in $S^1\times S^1$ with the additional property that the
map $(a,b,c,d)\in S^1\to \mu[a,b)\times [c,d)$ is
$\lambda$-measurable.

\section{Cross-ratios and bounded cohomology}

Let $X$ be an arbitrary set
with infinitely many
points. We define an
\emph{anti-symmetric cross ratio} for $X$ to
be a bounded function $[\,]$
on the space of quadruples of
pairwise distinct points in $X$ with the following additional
properties.
\begin{enumerate}
\item[i)] $[\xi,\xi^\prime,\eta,\eta^\prime]=
-[\xi^\prime,\xi,\eta,\eta^\prime]$.
\item[ii)] $[\eta,\eta^\prime,\xi,\xi^\prime]=
- [\xi,\xi^\prime,\eta,\eta^\prime]$.
\item[iii)] $[\xi,\xi^\prime,\eta,\eta^\prime]+
[\xi^\prime,\xi^{\prime\prime},\eta,\eta^\prime]=
[\xi,\xi^{\prime\prime},\eta,\eta^\prime]$.
\end{enumerate}
Note that $[\,]$ can naturally be extended by $0$ to the
set of quadruples of points in $X$ of the form
$(\xi,\xi,\eta,\zeta)$ where $\xi,\eta,\zeta$ are
pairwise distinct.

Let $\Gamma$ be an arbitrary countable group which
acts on $X$.
We define a
\emph{$\Gamma$-cross ratio} on $X$ to be an anti-symmetric cross
ratio on $X$ which is invariant under the action of $\Gamma$. Let
${\cal C}(\Gamma,X)$ be the vector space of all anti-symmetric
$\Gamma$-cross ratios on $X$, equipped with the supremums-norm
$\Vert \,\Vert$. Thus we have $\Vert
[\,]\Vert=\sup\{[\xi,\xi^\prime,\eta,\eta^\prime]\mid
\xi,\xi^\prime,\eta,\eta^\prime\in X$ are pairwise distinct$\}$,
and $({\cal C}(\Gamma,X),\Vert\,\Vert)$ is naturally a Banach space.

The second real bounded cohomology group $H_b^2(\Gamma,\mathbb{R})$ of
$\Gamma$ is a vector space of equivalence classes of bounded
functions on $\Gamma\times\Gamma\times \Gamma$ which are invariant
under the free left action of $\Gamma$. It can be
equipped with a natural pseudo-norm which assigns
to a cohomology class the infimum of the supremums-norms of any
representative of the class. For the second bounded cohomology,
this pseudo-norm turns out to be a norm
which provides $H_b^2(\Gamma,\mathbb{R})$ with
the structure of a Banach space. We have.

\bigskip

{\bf Proposition 2.1:} {\it There is a continuous linear map
$\Omega:{\cal C}(\Gamma,X)\to H_b^2(\Gamma,\mathbb{R})$.}

{\it Proof:} Let $[\,]\in {\cal C}(\Gamma,X)$ and let
$(\xi,\zeta,\eta)$ be a triple of pairwise distinct points in $X$.
Choose an arbitrary point
$\nu\in X$ which is
distinct from any of
the points in the triple and define \begin{equation}
\phi(\xi,\eta,\zeta)=\frac{1}{2}\bigl([\xi,\zeta,\eta,\nu]+
[\zeta,\eta,\xi,\nu]+[\eta,\xi,\zeta,\nu]\bigr).\end{equation} We
claim first that this does not depend on the choice of $\nu$.
Namely,
let $\nu^\prime\in X-\{\xi,\eta,\zeta\}$ be
arbitrary. We compute
\begin{align}
[\xi,\zeta,\eta,\nu]+ [\zeta,\eta,\xi,\nu]+ &
[\eta,\xi,\zeta,\nu]- [\xi,\zeta,\eta,\nu^\prime] -
[\zeta,\eta,\xi,\nu^\prime]- [\eta,\xi,\zeta,\nu^\prime]\\
=[\xi,\zeta,\nu^\prime,\eta] +
[\xi,\zeta,\eta,\nu]+ & [\zeta,\eta,\nu^\prime,\xi]+
[\zeta,\eta,\xi,\nu]+
[\eta,\xi,\nu^\prime,\zeta]+
[\eta,\xi,\zeta,\nu] \notag \\
=[\xi,\zeta,\nu^\prime,\nu] +&
[\zeta,\eta,\nu^\prime,\nu]+[\eta,\xi,\nu^\prime,\nu]=0 \notag
\end{align}
by the defining properties of an anti-symmetric cross ratio. This
shows our claim. As a consequence, the function $\phi$ is well
defined, moreover it is invariant under the action of the
group $\Gamma$. Its supremums norm is bounded in absolute value by
$3\Vert [\,]\Vert/2$. We extend $\phi$ by $0$ to
the set of all triples of points in $X$ for which
at least two of the points coincide.
Then $\phi$ induces a bounded function $\mu$
on $\Gamma^3$
which is invariant under the free left action of $\Gamma$ as
follows. Choose any point $\xi\in X$
and define $\mu(\gamma_1,\gamma_2,\gamma_3)=
\phi(\gamma_1\xi,\gamma_2\xi,\gamma_3\xi)$.

We next show that this function $\mu$ is a
cocycle, i.e. that it is alternating and
its image under the natural
coboundary map $L^\infty(\Gamma^3,\mathbb{R})
\to L^\infty(\Gamma^4,\mathbb{R})$ vanishes.
For this it is enough
to show the corresponding properties for the function
$\phi$. Namely, observe that
for a triple $(a,b,c)$ of pairwise distinct points in $X$
and any $d\in X-\{a,b,c\}$ we have
\begin{align}\label{triplesum}
2(\phi(a,c,b) &+\phi(a,b,c))=[a,b,c,d]+
[b,c,a,d]+[c,a,b,d]\\+[a,c,b,d]+&[c,b,a,d]+[b,a,c,d]=0=
2(\phi(a,c,b)+\phi(c,a,b)).\notag
\end{align}
which shows that $\phi$ is alternating.
Next we have to show that
\begin{equation}\label{triple}
\phi(b,c,d)-\phi(a,c,d)+\phi(a,b,d)-\phi(a,b,c)=0
\end{equation} for any
quadruple $(a,b,c,d)$ of points in $X$.
By our observation (\ref{triplesum}),
equation (\ref{triple}) is valid for every quadruple
$(a,b,c,d)$ of points in $X$ for which at least $2$ of the
points coincide. On the other hand, for every
quadruple $(a,b,c,d)$ of pairwise distinct points in $X$
we have
\begin{align}
2(\phi(a,b,c)-\phi(b,c,d))&  =
[a,c,b,d]+ [c,b,a,d]+[b,a,c,d]\\
-[b,d,c,a]-[d,c,b,a]- & [c,b,d,a] = 2[c,b,a,d]\notag
\end{align}
by the defining properties of an anti-symmetric cross ratio.
Together with equation (\ref{triplesum}) we deduce with
the same calculation as before that
$\phi(a,c,d)-\phi(a,b,d)=-\phi(c,a,d)+\phi(a,d,b)=
-[d,a,c,b]=-[c,b,a,d]$ and consequently
$\phi$ satisfies the cocycle equation (\ref{triple}).
Therefore the function $\mu$ on $\Gamma$
defines a bounded cocycle and hence a
second bounded
cohomology class $\Omega([\,])\in H_b^2(\Gamma,
\mathbb{R})$.

We next observe that the cohomology class $\Omega([\,])$ does not
depend on the choice of the point $\xi\in X$ which we used to
define the function $\phi$. Namely, let $\eta$ be a different point
in $X$ and define
$\mu^\prime(\gamma_1,\gamma_2,\gamma_3)=
\phi(\gamma_1\eta,\gamma_2\eta,\gamma_3\eta)$. For $\gamma\in
\Gamma$ write $\nu(\gamma)=[\gamma\xi,\eta,\xi,\gamma\eta]$. By
equation (\ref{triplesum}),(6) we have $\nu(\gamma)= \phi(\xi,\eta,
\gamma\xi)+\phi(\gamma\xi,\eta,\gamma\eta)$ and consequently
$\mu(e,\gamma_1,\gamma_2)- \mu^\prime(e,\gamma_1,\gamma_2)-
\nu(\gamma_1)-\nu(\gamma_1^{-1}\gamma_2)+ \nu(\gamma_2)$ is the
value of $\phi$ on the boundary of a singular polyhedron of
dimension 3 whose vertices consist of the points
$\xi,\eta,\gamma_1\xi,\gamma_1\eta,\gamma_2\xi, \gamma_2\eta$ and
whose sides contain the simplices with vertices
$\xi,\gamma_1\xi,\gamma_2\xi$ and $\eta,\gamma_1\eta,
\gamma_2\eta$ as well as three quadrangles with the same set of
vertices.
Since $\phi$ is a cocycle, the evaluation of $\phi$ on this
polyhedron vanishes and therefore the function $\mu-\mu^\prime$ is
the image of a $\Gamma$-invariant bounded function on $\Gamma^2$
under the coboundary map $L^\infty(\Gamma^2)\to
L^\infty(\Gamma^3)$. As a consequence, for every $[\,]\in {\cal
C}(\Gamma,X)$ the bounded cohomology class $\Omega([\,])\in
H_b^2(\Gamma,\mathbb{R})$ only depends on the cross ratio $[\,]$
but not on the choice of a point $\xi\in X$. The resulting map
$\Omega:{\cal C}(\Gamma,X)\to H_b^2(\Gamma,\mathbb{R})$ is clearly
linear and continuous with respect to the supremums norm on ${\cal
C}(\Gamma,X)$ and the Gromov norm on $H_b^2(\Gamma,\mathbb{R})$.
\qed

\bigskip

For a topological space $X$ define a \emph{continuous}
anti-symmetric cross ratio to be a continuous function
$[\,]$ on the space of quadruples of pairwise distinct points
in $X$ which is an anti-symmetric cross ratio in the
sense of our definition.
If $G$ is any locally compact $\sigma$-compact
topological group which acts on a locally compact
$\sigma$-compact topological space
$X$ as a group of homeomorphisms, then we define
a continuous anti-symmetric $G$-cross-ratio to be
a continuous anti-symmetric cross ratio on $X$ which
is invariant under the diagonal action of $G$.

Every locally compact $\sigma$-compact group $G$ acts freely on
itself by left translations. Thus we can use Proposition 2.1 for
$X=G$ to deduce that there is a continuous linear map from the
vector space of all continuous anti-symmetric $G$-cross ratios on
$G$ into $H_{cb}^2(G,\mathbb{R})$. Namely, we have.

\bigskip

{\bf Corollary 2.2:} {\it Let $G$ be a locally compact
$\sigma$-compact group; then
$H_{cb}^2(G,\mathbb{R})$ is the quotient of
the space of continuous
anti-symmetric $G$-cross ratios on $G$
under the coboundary
relation.}

{\it Proof:} By Proposition 2.1 and its proof there is a
linear map $\Omega$ from the space of all
continuous anti-symmetric $G$-cross ratios on $G$
into $H_{cb}^2(G,\mathbb{R})$. We have to show that
$\Omega$ is surjective.
For this recall first
that every second continuous
bounded cohomology class for $G$ is the
class of a $G$-invariant continuous \emph{alternating} function
$\phi$ on $G^3$ which
satisfies the cocycle equation (\ref{triple})
(see \cite{I,M}).
For such a function $\phi$ define
$[g,g^\prime,h,h^\prime]=\phi(h,g^\prime,g)-\phi(g^\prime,g,h^\prime)$
as in the proof of Proposition 2.1. Then $[\,]$ is continuous
and $G$-invariant;
we claim that it satisfies the requirements of
an anti-symmetric $G$-cross ratio.

Namely, first of all we have
\begin{equation}
[b,c,a,d]=\phi(a,c,b)-\phi(c,b,d)=-\phi(a,b,c)+\phi(b,c,d)=-[c,b,a,d]
\end{equation}
since $\phi$ is alternating
and similarly
\begin{equation}
[a,d,c,b]=\phi(c,d,a)-\phi(d,a,b)=-\phi(a,b,d)+\phi(a,c,d)=-[c,b,a,d].
\end{equation}
In the same way we obtain from the cocycle equation
(\ref{triple}) for $ \phi$ that
\begin{equation}
[a,b,c,d] +[b,b^\prime,c,d]=[a,b^\prime,c,d]
\end{equation}
which shows that $[\,]$ is indeed a $G$-invariant
anti-symmetric cross ratio. By the proof of
Proposition 2.1, the
cohomology class of $\phi$ is just $\Omega([\,])$.
Thus the map $\Omega$ is surjective.
\qed

\bigskip

A locally compact
$\sigma$-compact topological group $G$
admits a \emph{strong boundary} \cite{K03} which
is a standard measure space $(B,\mu)$ together
with a measure class preserving action of $G$ such
that the following two conditions are satisfied.
\begin{enumerate}
\item The $G$-action on $B$ is amenable.
\item For any separable Banach-$G$-module $E$, any
measurable $G$-equivariant map $B\times B\to E$ is
essentially constant.
\end{enumerate}
Define a \emph{measurable anti-symmetric cross ratio}
on $B$ to be a \emph{measurable} essentially bounded
$G$-invariant function on
the space of quadruples of distinct points in $B$ which
satisfies almost everywhere the defining equation for a cross ratio.
Our above considerations immediately imply.

\bigskip

{\bf Lemma 2.3:} {\it Let $(B,\mu)$ be a strong boundary for $G$;
then $H_{cb}^2(G,\mathbb{R})$ is naturally isomorphic to
the space of $G$-invariant
$\mu$-measurable anti-symmetric cross ratios
on $B$.}

{\it Proof:} Let $G$ be a locally compact $\sigma$-compact
group and let $(B,\mu)$ be a strong boundary
for $G$.
Denote for $k\geq 1$ by
$L^\infty(B^k,\mu^k)$ the Banach space of essentially
bounded $\mu^k$-measurable functions on $B^k$; it
contains the closed subspace $L^\infty(B^k,\mu^k)^G$
of all $G$-invariant such functions. By the
results of Burger and Monod \cite{BM02},
$H_{cb}^2(G,\mathbb{R})$ coincides with the
second bounded cohomology group of the resolution
\begin{equation}
\mathbb{R}\to L^\infty(B,\mu)^G\to
L^\infty(B^2,\mu^2)^G\to \dots
\end{equation}
with the usual homogeneous coboundary operator \cite{M}.
By Proposition 2.1 and its proof, there is a natural bijection
between the space of
$\mu^4$-measurable essentially bounded
anti-symmetric $G$-cross ratios on $B$ and the space of
$\mu^3$-measurable bounded cocycles on $B$, i.e.
the space of alternating functions in
$L^\infty(B^3,\mu^3)^G$ which satisfy
the cocycle equation. Any two such bounded cocycles are
cohomologous if and only if they differ by the image under
the coboundary operator of an element of $L^\infty(B^2,
\mu^2)^G$. However, by ergodicity of the measure class
$\mu^2$ under the diagonal action of $\Gamma$, the
vector space $L^\infty(B^2,\mu^2)^G$ only contains
essentially
constant functions. This means that the natural map which
assigns to a $\mu^4$-measurable $G$-cross ratio its corresponding
bounded cohomology class is injective.
\qed

\bigskip

We next look at a specific example. Namely,
the ideal boundary of the oriented hyperbolic plane ${\bf H}^2$
can naturally be identified with the oriented circle $S^1$.
The action on $S^1$ of the group $PSL(2,\mathbb{R})$
of orientation preserving isometries of ${\bf H}^2$
preserves the measure class of the Lebesgue measure
$\lambda$. Moreover, $(S^1,\lambda)$ is a strong
boundary for $PSL(2,\mathbb{R})$.
A triple $(a,b,c)$ of pairwise distinct points in $S^1$ will
be called \emph{ordered} if the point $b$ is contained in the
interior of the oriented subinterval $(a,c)$ of $S^1$.
Similarly
we define a quadruple $(a,b,c,d)$ or pairwise distinct points
on $S^1$ to be ordered.
Ordered quadruples define a $PSL(2,\mathbb{R})$-invariant
subset of the space of quadruples of pairwise distinct points in $S^1$.
There is a unique antisymmetric cross ratio $[\,]_0$ on
$S^1$ with the following properties.
\begin{enumerate}
\item $[a,b,c,d]_0=0$ if $(a,b,c,d)$ is ordered.
\item $[a,b,c,d]_0=1$ if $(a,c,b,d)$ is ordered.
\end{enumerate}
Clearly this cross ratio is measurable and invariant under
the action of $PSL(2,\mathbb{R})$.
We have.

{\bf Lemma 2.4:} {\it \begin{enumerate}
\item Up to a constant, $[\,]_0$ is the
unique non-trivial antisymmetric cross ratio on $S^1$ which
is invariant under the full group $PSL(2,\mathbb{R})$ of
orientation preserving isometries of ${\bf H}^2$.
\item If $\mu$ is any Borel measure on $S^1$ without atoms
and if
$[\,]$ is any $\mu$-measurable
antisymmetric cross ratio on $S^1$ which satisfies
$[a,b,c,d]=0$ for $\mu$-almost every ordered quadruple of
pairwise distinct points in $S^1$ then $[\,]$ is a multiple
of $[\,]_0$.
\end{enumerate}
}

{\it Proof:} If $[\,]$ is any antisymmetric cross ratio which is
invariant under the full group $PSL(2,\mathbb{R})$ then the same
is true for the induced function $\phi$ on the set of triples of
pairwise distinct points in $S^1$. Recall that the group
$PSL(2,\mathbb{R})$ acts transitively on the space of
\emph{ordered} triples of distinct points in $S^1$ and
consequently we have $\phi(a,b,c)=\phi(b,c,d)$ for every ordered
quadruple $(a,b,c,d)$
of pairwise distinct points in $S^1$. Hence equation (6)
in the proof of Proposition 2.1 implies that $[a,b,c,d]=0$
whenever the quadruple $(a,b,c,d)$ of pairwise distinct points in
$S^1$ is ordered. As a consequence, if
$[\,]\not\equiv 0$ then after possibly multiplying $[\,]$ with a
constant we conclude from the definition of an anti-symmetric
cross ratio that $(a,b,c,d)=1$ whenever $(a,c,b,d)$ is ordered.
In other words, up to a constant
we have $[\,]=[\,]_0$ which shows the first part
of the lemma.

To show the second part of the lemma, let $\mu$ be a Borel measure
on $S^1$ without atoms and let $[\,]$ be any $\mu$-measurable
antisymmetric cross ratio on $S^1$ such that $[a,b,c,d]=0$ for
almost every ordered quadruple of pairwise distinct points in
$S^1$. Let $(a,b,c,d)$ be a quadruple of pairwise distinct points
in $S^1$ such that $[a,b,c,d]\not=0$. By the symmetry relation of
an antisymmetric cross ratio we may assume that $(a,c,b,d)$ is
ordered. Thus the point $b$ is contained in the oriented open
subinterval $(c,d)$ of $S^1$. Now if $b^\prime\in (c,d)$ is
another such point then by the defining properties of an
antisymmetric cross ratio we have
$[a,b,c,d]+[b,b^\prime,c,d]=[a,b^\prime,c,d]$. Assume first that
$b^\prime\in (b,d)$; then the quadruple $(b,b^\prime,d,c)$ is
ordered and consequently we have $[b,b^\prime,c,d]=0$ by
assumption and $[a,b,c,d]=[a,b^\prime,c,d]$. Via exchanging $b$
and $b^\prime$ we conclude that $[a,b,c,d]$ does not depend on
$b\in (c,d)$. Similarly we deduce that in fact $[a,b,c,d]$ only
depends on the order of the points in the quadruple $(a,b,c,d)$.
In other words, $[\,]$ is a multiple of $[\,]_0$. \qed
\bigskip

Lemma 2.4 shows in particular that $[\,]_0$ corresponds to a
generator $\alpha$ of the group $H_{cb}^2(PSL(2,\mathbb{R}),\mathbb{R})$.

Now let $\Gamma<PSL(2,\mathbb{R})$ be any finitely generated
torsion free discrete group. Then $\Gamma$ acts freely on the
hyperbolic plane ${\bf H}^2$ as a properly discontinuous group of
isometries, and $S={\bf H}^2/\Gamma$ is an oriented surface.
The unit circle $S^1=\partial {\bf H}^2$ with the
$\Gamma$-invariant measure class of the Lebesgue measure $\lambda$
is a strong boundary for $\Gamma$. If $\Gamma$ is cocompact then
we have $H^2(\Gamma,\mathbb{R})=H^2(M,\mathbb{R})=\mathbb{R}$ and
there is a natural surjective restriction homomorphism
$H^2_b(\Gamma,\mathbb{R})\to H^2(\Gamma,\mathbb{R})$. The diagram
\[\begin{CD}
H_b^2(\Gamma,\mathbb{R}) @>>>  H^2(\Gamma,\mathbb{R})\\                         @AAA                            @AAA \\
H^2_{cb}(PSL(2,\mathbb{R}),\mathbb{R}) @>>>
H_c^2(PSL(2,\mathbb{R}),\mathbb{R})
\end{CD}
\]
commutes and hence in our identification of
$H^2_b(\Gamma,\mathbb{R})$ with the space of $\Gamma$-invariant
antisymmetric cross ratios on the unit circle $S^1$, the
$PSL(2,\mathbb{R})$-invariant cross ratio $[\,]_0$ from Lemma 2.4
maps to a multiple of the Euler class of the tangent bundle of
$S$. An easy calculation shows that the multiplication factor is
in fact 1 (compare \cite{BI05} for a detailed explanation of this
fact). If $\Gamma$ is not cocompact then $S$ is homeomorphic to
the interior of a compact oriented surface $\hat S$ whose boundary
$\partial \hat S$ is a finite union of circles. In particular, the
fundamental group of each of the boundary circles is amenable. Now
the bounded cohomology of any countable cellular space $X$ is
naturally isomorphic to the bounded cohomology of its fundamental
group $\pi_1(X)$ \cite{I}. Thus $H_b(\partial \hat
S,\mathbb{R})=\{0\}$ and hence the relative group $H^2_b(\hat
S,\partial \hat S,\mathbb{R})$ is isomorphic to $H_b^2(\hat
S,\mathbb{R})= H^2_b(\Gamma,\mathbb{R})$. The image under the
natural map $\iota:H^2_b(\hat S,\partial \hat S,\mathbb{R})\to
H^2(\hat S,\partial \hat S,\mathbb{R})$ of the pull-back
$i^*\alpha$ of $\alpha$ under the inclusion $i:\Gamma\to
PSL(2,\mathbb{R})$ is just the Euler class of the (bordered)
surface $\hat S$. In other words, the evaluation of
$\iota(i^*\alpha)$ on the fundamental class of the pair $(\hat
S,\partial \hat S)$ just equals the Euler characteristic of $S$
(we refer to \cite{BI05} for a detailed discussion of a slightly
different viewpoint on these facts).

To describe the kernel $Q$ of the natural map
$H_b^2(\Gamma,\mathbb{R})\to
H^2(\hat S,\partial \hat S,\mathbb{R})$, let $\Delta$
be the diagonal in $S^1\times S^1$ and recall that a finitely
additive signed measure on $S^1\times S^1-\Delta$ is a function
$\mu$ which assigns to each oriented rectangle of the form
$[a,b)\times [c,d)$ where $a,b,c,d$ are pairwise distinct a number
$\mu[a,b)\times [c,d)\in \mathbb{R}$ such that
\begin{equation}
\mu[a,b)\times
[c,d)+\mu[b,b^\prime)\times [c,d)= \mu[a,b^\prime)\times [c,d)
\end{equation}
whenever $[a,b^\prime)=[a,b)\cup [b,b^\prime)$ is
disjoint from $[c,d)$.
This measure is flip-anti\-in\-va\-riant if $\mu[a,b)\times
[c,d)=-\mu[c,d)\times [a,b)$. Let $\lambda$ be the Lebesgue measure
on $S^1$ which defines the unique
$PSL(2,\mathbb{R})$-invariant measure class.
The collection of all
finitely additive flip-antiinvariant signed
measures $\mu$ on $S^1\times S^1$ with the additional property
that the function $(a,b,c,d)\to \mu[a,b)\times [c,d)$ is
$\lambda$-measurable is clearly a vector space; we call it the
space of \emph{$\lambda$-measurable} flip-antiinvariant
signed measures. We have (compare
Section 2 of \cite{H05a}).

\bigskip

{\bf Corollary 2.5:} {\it There is a linear isomorphism $\Psi$ of
the kernel of the map $H_b^2(\Gamma,\mathbb{R})\to H^2(\hat
S,\partial \hat S,\mathbb{R})$ onto the vector space of all
$\lambda$-measurable finite flip-antiinvariant $\Gamma$-invariant
finitely additive signed measures on $S^1\times S^1-\Delta$.}

{\it Proof:} Every $\Gamma$-invariant $\lambda$-measurable
anti-symmetric cross ratio $[\,]$ on $S^1$ defines a finite
finitely additive measurable signed measure $\mu=\Psi([\,])$ by
assigning to an ordered quadruple $(a,b,c,d)$ of pairwise distinct
points in $S^1$ the value $\mu[a,b)\times [c,d)=[a,b,c,d]$. The
map $\Psi:[\,]\to\Psi([\,])$ is clearly linear, moreover by Lemma
2.4 its kernel is spanned by the cross ratio $[\,]_0$ which
defines the Euler class of $S$ viewed as a class in $H^2(\hat
S,\partial \hat S,\mathbb{R})$. As as consequence, the restriction
of the map $\Psi$ to the kernel $Q$ of the natural map
$H_b^2(\Gamma,\mathbb{R})\to H^2(\hat S,\partial \hat
S,\mathbb{R})$ is injective.

Now let $\mu$ be any $\lambda$-measurable finitely additive
$\Gamma$-invariant flip anti-invariant signed measure on
$S^1\times S^1-\Delta$. Define an antisymmetric cross ratio $[\,]$
as follows. First, if $(a,b,c,d)$ is ordered and if $a,b,c,d$ are
typical points for $\lambda$ such that $\mu[a,b)\times [c,d)$ is
defined then we define $[a,b,c,d]= \mu[a,b)\times [c,d)$. This
determines $[a,b,c,d]$ whenever $(a,b,c,d)$ or $(b,a,c,d)$ or
$(a,b,d,c)$ is ordered. Choose an arbitrary fixed typical ordered
quadruple $(a,b,c,d)$ and define $[a,c,b,d]=0$. As in the proof of
Lemma 2.4, if $b^\prime\in (a,b)$ is arbitrary (and typical) then
$(c,a,b^\prime,b)$ is ordered and we define
$[a,c,b^\prime,d]=[a,c,b^\prime,b]+[a,c,b,d]= -\mu [c,a)\times
[b^\prime,b)$. Similarly, for a point $b^\prime\in (a,b)$ and
every $a^\prime\in (d,a)$ define $[a^\prime,c,b^\prime,
d]=[a^\prime,a,b^\prime,d]+[a,c,b^\prime,d]=\mu[a^\prime,a)\times
[b^\prime,d)+[a,c,b^\prime,d]$. Successively we can define in this
way an anti-symmetric measurable $\Gamma$-cross ratio on $S^1$. As
a consequence, the map $\Psi$ is surjective. This shows the
corollary. \qed

\bigskip

It follows from
our consideration that there is a direct
decomposition $H_b^2(\Gamma,\mathbb{R})=
Q\oplus {\rm ker}(\Psi)$ where
$Q$ is the kernel of the natural map $H_b^2(\Gamma,\mathbb{R})
\to H^2(\hat S,\partial \hat S,\mathbb{R})$ and
$\Psi$ is the map which
associates to a $\lambda$-measurable anti-symmetric
$\Gamma$-cross ratio $[\,]$ its corresponding
finitely additive signed measure on $S^1\times S^1-\Delta$.

\bigskip

{\bf Remark:} Let $\Gamma$ be a cocompact torsion free lattice in
$PSL(2,\mathbb{R})$. Then $S={\bf H}^2/\Gamma$ is a closed surface
of genus $g\geq 2$. The \emph{geodesic flow} $\Phi^t$ acts on the
unit tangent bundle $T^1S$ of $S$. Call two H\"older continuous
functions $f,u:T^1S\to \mathbb{R}$ \emph{cohomologous} if the
integrals of $f,u$ over every periodic orbit for $\Phi^t$
coincide. This defines an equivalence relation on the vector space
of all H\"older function on $T^1S$. The \emph{flip} $v\to -v$ acts
on $T^1M$ and induces a map on the space of cohomology classes. We
showed in \cite{H05a} that the space of cohomology classes of
\emph{flip anti-invariant} H\"older functions naturally embeds
into $H_b^2(\Gamma,\mathbb{R})$. Our above consideration yields
another such embedding. Namely, every such function $f$ admits a
unique \emph{Gibbs equilibrium state} which defines a
$\Gamma$-invariant locally finite measure $\mu$ on the
complement of the diagonal in $S^1\times S^1$. The function
$(a,b,c,d)\to \mu[a,b)\times [c,d)$ is H\"older continuous. Since
$f$ is anti-invariant by assumption, the measure is \emph{not}
invariant under the flip unless $f$ is cohomologous to $0$. In
particular, the signed measure $\tilde \mu$ defined by $\tilde
\mu[a,b)\times [c,d)=\mu[a,b)\times [c,d)-\mu[c,d)\times [a,b)$ is
non-trivial, flip anti-invariant and H\"older continuous and hence
by Corollary 2.5 it defines a non-trivial class in the kernel of
the natural map $H_b^2(\Gamma,\mathbb{R})\to
H^2(\Gamma,\mathbb{R})$. We do not know whether the signed
measures of this form are dense in the space of all finitely
additive $\lambda$-measurable signed measures equipped with the
supremums norm.

\section{Large groups of isometries of hyperbolic spaces}

In this section we consider a proper hyperbolic geodesic metric
space $X$. The \emph{Gromov boundary} $\partial X$ of $X$ is
defined as follows. For a fixed point $x\in X$, define the
\emph{Gromov product} $(y,z)_{x}$ based at $x$ of two points
$y,z\in X$  by
\begin{equation}
(y,z)_{x}=\frac{1}{2}\bigl(d(y,x)+d(z,x)-d(y,z)\bigr).
\end{equation}
Call two sequences $(y_i),(z_j)\subset X$ \emph{equivalent} if
$(y_i,z_i)_x\to \infty$ $(i\to \infty)$. By hyperbolicity of $X$,
this notion of equivalence defines an equivalence relation for the
collection of all sequences $(y_i)\subset X$ with the additional
property that $(y_i,y_j)_{x}\to \infty$ $(i,j\to\infty)$ \cite{BH}.
The boundary $\partial X$ of $X$ is the set of equivalence classes
of this relation.

There is a natural topology on $X\cup \partial X$ which restricts
to the given topology on $X$.
With respect to this topology, a
sequence $(y_i)\subset X$ converges to $\xi\in
\partial X$ if and only if we have $(y_i,y_j)_{x}\to \infty$ and the
equivalence class of $(y_i)$ defines $\xi$. Since $X$ is proper
by assumption, the space
$X\cup \partial X$ is compact and metrizable.
Every isometry of $X$ acts
naturally on $X\cup
\partial X$ as a homeomorphism.
Moreover, for every $x\in X$ and every $a\in \partial X$ there
is a geodesic ray $\gamma:[0,\infty)\to X$ with $\gamma(0)=x$
and $\lim_{t\to\infty}\gamma(t)=a$.

For every proper metric space $X$ the isometry group ${\rm
Iso}(X)$ of $X$ can be equipped with a natural locally compact
$\sigma$-compact metrizable topology, the so-called \emph{compact
open topology}. With respect to this topology, a sequence
$(g_i)\subset {\rm Iso}(X)$ converges to some isometry $g$ if and
only if $g_i\to g$ uniformly on compact subsets of $X$.
In this topology, a
closed subset $Y\subset {\rm Iso}(X)$ is compact if and only if
there is a compact subset $K$ of $X$ such that $gK\cap
K\not=\emptyset$ for all $g\in Y$. In particular, the action of
${\rm Iso}(X)$ on $X$ is proper.
In the
sequel we always assume that subgroups of ${\rm Iso}(X)$ are
equipped with the compact open topology.

The \emph{limit set} $\Lambda$ of a subgroup $G$ of the isometry
group of a proper hyperbolic geodesic metric space $X$ is the set
of accumulation points in $\partial X$ of one (and hence every)
orbit of the action of $G$ on $X$. If $G$ is non-compact then its
limit set is a compact non-empty $G$-invariant subset of $\partial
X$. The group $G$ is called \emph{elementary} if its limit set
consists of at most two points. In particular, every compact
subgroup of ${\rm Iso}(X)$ is elementary. If $G$ is non-elementary
then its limit set $\Lambda$ is uncountable without isolated
points, and $G$ acts as a group of homeomorphisms on $\Lambda$.
An element $g\in
{\rm Iso}(X)$ is called \emph{hyperbolic} if its action on
$\partial X$ admits an attracting fixed point $a\in \Lambda$ and a
repelling fixed point $b\in\Lambda-\{a\}$ and if it acts on $\Lambda$
with \emph{north-south-dynamics} with respect to these fixed
points.

Assume that $G<{\rm Iso}(X)$ is a closed non-elementary subgroup
of ${\rm Iso}(X)$ with limit set $\Lambda\subset\partial X$. If
the group $G$ does not act transitively on the complement of the
diagonal in $\Lambda\times\Lambda$ then the kernel of the natural
map $H_{cb}^2(G,\mathbb{R})\to H_c^2(G,\mathbb{R})$ is infinite
dimensional \cite{H05b}. Extending Lemma 6.1 of \cite{BM99}, we
prove Theorem B from the introduction and show that this is not
true for groups $G$ which act transitively on the complement of the
diagonal in $\Lambda\times \Lambda$.

\bigskip

{\bf Proposition 3.1:} {\it Let $G<{\rm Iso}(X)$ be a closed
non-elementary subgroup with limit set $\Lambda$. If $G$ acts
transitively on the complement of the diagonal in $\Lambda\times
\Lambda$ then the kernel of the natural map
$H_{cb}^2(G,\mathbb{R})\to H_c^2(G,\mathbb{R})$ is trivial.}

{\it Proof:} Let $G<{\rm Iso}(X)$ be a closed non-elementary
subgroup which acts transitively on the space $A$ of pairs of
distinct points in $\Lambda$. Then for every fixed point $a\in
\Lambda$ the stabilizer $G_a$ of $a$ acts transitively on
$\Lambda-\{a\}$.

Following \cite{BM99}, an element in the kernel
$H_{cb}^2(G,\mathbb{R})\to H_c^2(G,\mathbb{R})$ can
be represented by a \emph{continuous quasi-morphism},
i.e. a continuous function $q:G\to \mathbb{R}$ which satisfies
\begin{equation}
\sup_{g,h\in G}\vert q(g)+q(h)-q(gh)\vert <\infty.
\end{equation}
Such a continuous quasi-morphism is bounded on
every \emph{compact} subset of $G$ and is bounded
on each fixed conjugacy class in $G$. The quasi-morphism $q$ defines
a non-trivial element of $H_{cb}^2(G,\mathbb{R})$ only
if $q$ is unbounded.

Since $G$ is non-elementary by assumption, $G$ contains a
hyperbolic element $g\in G$. Let $a\in \Lambda$ be the attracting
fixed point of $g$ and let $b\in \Lambda-\{a\}$ be the repelling
fixed point. Then $g$ preserves the set of geodesics connecting
$b$ to $a$; note that such a geodesic exists since $X$ is proper.
For $x,y\in X$ and a geodesic ray
$\gamma:[0,\infty)\to X$ connecting $x$ to $a$ write
$\beta(y,\gamma)=\lim\sup_{t\to \infty}(d(y,\gamma(t))-t)$ and
define the \emph{Busemann function}
\begin{equation}
\beta_a(y,x)=\sup\{\beta(y,\gamma)\mid \gamma\,\text{is a
geodesic ray connecting}\,x\,\text{to}\,a\}.
\end{equation}
Then we have $\vert \beta_a(y,x)\vert \leq d(y,x)$,
moreover
there is a number $c>0$ only
depending on the hyperbolicity constant for $X$ such that
\begin{equation}\label{buse}
\vert \beta_a(\cdot,y)-\beta_a(\cdot,x)+\beta_a(y,x)\vert \leq c
\end{equation}
for all $x,y\in X$ (Proposition 8.2 of \cite{GH}).
Let $\gamma:[0,\infty)\to X$
be a geodesic ray connecting $x$ to $a$.
By hyperbolicity, for every $t>0$
every geodesic ray connecting $x$ to $a$ passes through
a neighborhood of $\gamma(t)$ of uniformly
bounded diameter and therefore via possibly enlarging $c$
we may assume that \begin{equation}\label{distance}
\vert \beta_a(\gamma(t),x)+d(\gamma(t),x)\vert \leq c
\end{equation}
independent of $x,\gamma,t$. Define the \emph{horosphere} at $a$
through $x$ to be the set $H_a(x)=\beta_a(\cdot,x)^{-1}(0)$. For
all $x,y\in X$ the distance between the horospheres
$H_a(x),H_a(y)$ is not smaller than $\vert \beta_a(x,y)\vert -c$.
Moreover, the estimates (\ref{buse}),(\ref{distance}) show that
for all $x,y$ with $\beta_a(x,y)\leq 0$ the distance between $y$
and $H_a(x)$ is not bigger than $\vert \beta_a(x,y)\vert +2c$.

Let $W(a,b)\subset X$
be the closed non-empty subset of all points in $X$ which lie on
a geodesic connecting $b$ to $a$.
The isometry $g\in G$ is hyperbolic with fixed points $a,b\in \Lambda$
and therefore it preserves $W(a,b)$.
If we denote by $\Gamma$ the infinite cyclic
subgroup of $G$ generated by $g$ then $W(a,b)/\Gamma$ is
compact.
As a consequence,
there is a number $\nu>0$ and for every
$x\in W(a,b)$
and every $t\in \mathbb{R}$ there is a number $k(t)\in
\mathbb{Z}$ with $\vert \beta_a(g^{k(t)}x,x)-t\vert <\nu$. It
follows from this and (\ref{distance})
that  for every $y\in X$ there is some $k=k(y)\in
\mathbb{Z}$ such that the distance between $g^kx$ and the horosphere
$H_a(y)$ is at most $\delta_0=\nu+2c$.

By assumption, the stabilizer $G_a<G$ of the point $a\in \Lambda$
acts transitively on $\Lambda-\{a\}$. Thus there is for every
$\zeta\in \Lambda-\{a\}$ an element $h_\zeta\in G_a$ with
$h_\zeta(b)=\zeta$. Let $x\in W(a,b)$;
since $h_\zeta\in G_a$ we
have $h_\zeta^{-1}(H_a(x)) =H_a(h_\zeta^{-1}(x))$. By our above
consideration, there is a number $\ell\in \mathbb{Z}$ such that
the distance between $g^{\ell}(x)$ and $H_a(h_\zeta^{-1}(x))$
is at most $\delta_0$ and therefore the
distance between $h_\zeta\circ g^\ell(x)$ and
$H_a(x)$ is at most $\delta_0$. Hence via
replacing $h_\zeta$ by $h_\zeta\circ g^\ell\in G_a$ we may assume
that the distance between $h_\zeta(x)$ and $H_a(x)$
is at most $\delta_0$. In particular, we have $\vert
\beta_a(\cdot,x)-\beta_a(\cdot, h_\zeta(x))\vert \leq \delta_0+c$.
Now $h$ maps the geodesic $\gamma$ connecting $x$ to $a$ to a
geodesic $h\gamma$ connecting $h_\zeta x$ to $a$ and hence by the
definition of the Busemann functions, for every $t\in \mathbb{R}$
the distance between $h_\zeta(\gamma(t))$ and
$H_a(\gamma(t))$ is bounded from above by $\delta_0+3c$.
Since $\Gamma$ preserves the set of geodesics connecting
$b$ to $a$ we conclude that there is a universal constant
$\delta_1>0$ such that for every $\ell\geq 0$ the distance
between $h_\zeta(g^\ell x)$ and $H_a(g^\ell(x))$ is at most
$\delta_1$.

For $x\in X$
denote by $N_{a,x}\subset G_a$ the set of all elements $h\in G_a$ with
the property that the distance between $H_a(x)$ and
$hx$ is at most $\delta_0$. Our above consideration
shows that for every $x\in W(a,b)$
and every $\zeta\in \Lambda-\{a\}$ there is
some $h_\zeta\in N_{a,x}$ which maps $b$ to $\zeta$.
For every $h\in N_{a,x}$ the sequence
$(g^{-\ell}\circ h\circ g^{\ell})_\ell$ is contained in $G_a$ and
by the above consideration, it maps $x$ into
the $\delta_1$-neighborhood of $H_a(x)$. Moreover, the sequence
$(g^{-\ell}\circ h\circ g^\ell(b))_\ell\subset \Lambda$ converges
as $\ell\to \infty$ to $b$. This means that there is a number
$\delta_2>\delta_1$ such that
for sufficiently large
$\ell$ the element $g^{-\ell}\circ h\circ g^{\ell}\in G_a$ maps
the point $x$ into the closed $\delta_2$-neighborhood $B_x$ of $x$.
The group
$\Gamma$ acts on $W(a,b)$ cocompactly and hence if
$C\subset W(a,b)$ is a compact
fundamental domain for this action then
$B=\cup_{x\in C}B_x$ is compact, moreover for every $x\in W(a,b)$,
every element $h\in N_{a,x}$
is conjugate in
$G$ to an element in the compact subset $K=\{u\in G\mid
uB\cap B\not=\emptyset\}$ of $G$.
As a consequence, the restriction to
$N_{a,b}=\cup_{x\in  W(a,b)}N_{a,x}$ of any continuous
quasi-morphism $q$ on $G$ is uniformly bounded.
By our
assumption on $G$ the sets $N_{a,b}$ $((a,b)\in A)$ are pairwise
conjugate in $G$ and hence $q$ is uniformly bounded on
$\cup_{(a,b)\in A}N_{a,b}=N$.

Next we show that the restriction of a quasi-morphism $q$ to the
subgroup $G_{a,b}$ of $G_a$ which stabilizes both points $a,b\in
\Lambda$ is bounded. For this consider an arbitrary element $u\in
G_{a,b}$. Let $x\in W(a,b)$ and let $\beta_a(ux,x)=\tau$; we may
assume that $u\not\in N_{a,b}$ and hence after possibly exchanging
$b$ and $a$ that $\tau<0$. 
Since $u$ preserves $W(a,b)$, by 
(\ref{distance}) we have $\vert d(ux,x)+\tau\vert \leq \tilde c$
where $\tilde c>0$ is a universal constant. Choose some $h\in G$
with $h(a,b)=(b,a)$; then $h$ preserves the set $W(a,b)$. Let
$g\in G_{a,b}$ be the fixed hyperbolic element as before. By our
above consideration, via composing $h$ with $g^\ell$ for a number
$\ell\in \mathbb{Z}$ depending on $h$ we may assume that
$d(x,hx)\leq \delta_2$ for a universal constant $\delta_2>0$. We
have $hux\in W(a,b)$ and $\vert \beta_b(hux,x)-\beta_a(ux,x)\vert\leq
\vert\beta_b(hux,hx)-\beta_a(ux,x)\vert
+d(hx,x)\leq \delta_2$. On the other
hand, from another application of the estimate (\ref{distance}) we
obtain the existence of a constant $\delta_3>0$ such that
$\vert\beta_b(z,y)-\beta_a(y,z)\vert \leq \delta_3$ for all
$y,z\in W(a,b)$. But this just means that the distance between
$hux$ and $u^{-1}x$ is uniformly bounded and hence the
distance between $uhux$ and $x$ is uniformly bounded as
well. As a consequence, the element $uhu$ is contained in a
fixed compact subset of $G$. As before, we conclude from this that
$\vert q(u)\vert $ is bounded from above by a universal constant
not depending on $u$ and hence the restriction of $q$ to $G_{a,b}$
is uniformly bounded.

Now let $h\in G$ be arbitrary and assume that $ha=x,hb=y$. We
showed above that there are $h_y\in N\cap G_x,h_x\in N\cap G_b$
with $h_y(y)=b$ and $h_x(x)=a$; then $h^\prime=h_xh_yh\in G_{a,b}$
and $\vert q(h^\prime)-q(h)\vert $ is uniformly bounded. As a
consequence, $q$ is bounded and hence it defines the trivial
bounded cohomology class. This completes the proof of the
proposition. \qed

\bigskip

From now on we assume that $X$ is of \emph{bounded growth} which
means that there is a number $b>1$ such that for every $R>1$,
every metric ball of radius $R$ contains at most $be^{bR}$
disjoint metric balls of radius $1$. Let $G<{\rm Iso}(X)$ be a
closed non-elementary group with limit set $\Lambda$.
Let $(B,\mu_0)$ be a strong boundary for $G$. By Lemma 2.2 of
\cite{H05b} there is a $G$-equivariant map
$\phi:B\to \Lambda$. The image $\mu$ of $\mu_0$ under the map $\phi$
is a measure on $\Lambda$ whose measure class is
invariant under the action of $G$ and such that the diagonal
action of $G$ on $\Lambda\times\Lambda$ is ergodic with respect to
the class of $\mu\times\mu$. Since
$X$ is of bounded growth by assumption, the stabilizer
in ${\rm Iso}(X)$ of each
point in $\partial X$ is amenable and therefore the measure space
$(\Lambda,\mu)$ is a strong boundary for $G$ (compare the
discussion in \cite{A94,A96,K03}). In particular, Corollary 2.3
implies that $H_{cb}^2(G,\mathbb{R})$ is isomorphic to the space
of $\mu$-measurable anti-symmetric $G$-cross ratios on $\Lambda$.

Now assume that $G$
acts transitively on the complement of the diagonal
in $\Lambda\times \Lambda$. The stabilizer
$G_{a,b}$ of $(a,b)$ acts on $\Lambda-\{a,b\}$
as a group of homeomorphisms. Moreover,
there is some $g\in G$ which
maps $(a,b)$ to $(b,a)$.
Define the group $G$ to be \emph{directed} if there is a
$G_{a,b}$-invariant subset $U(a,b)$ of $\Lambda-\{a,b\}$ of
positive $\mu$-mass such that $gU(a,b)\cap U(a,b)=\emptyset$. We
have.

\bigskip

{\bf Lemma 3.2:} {\it Let $G<{\rm Iso}(X)$ be a closed
subgroup with limit set $\Lambda$. Assume that
$G$ acts transitively on the complement $A$ of the diagonal
in $\Lambda\times\Lambda$; if $H_{cb}^2(G,\mathbb{R})\not=\{0\}$
then $G$ is directed.}

{\it Proof:} Let $G<{\rm Iso}(X)$ be a closed subgroup with limit
set $\Lambda$ which acts transitively on the complement
$A$ of the diagonal in $\Lambda\times \Lambda$. If
$H_{cb}^2(G,\mathbb{R})\not=0$ then there is some $G$-invariant
measurable nontrivial
alternating function $\phi$ on the space of triples of
pairwise distinct points in $\Lambda$ which satisfies the cocycle
identity (\ref{triplesum}).

Let $\phi$ be any $G$-invariant $\mu$-measurable alternating
function on the space of triples of pairwise distinct points in
$\Lambda$.
Let $(a,b)\in A$; since the action of $G$ on $A$ is transitive,
there is a set $U\subset \Lambda$ of positive $\mu$-mass such that
$\phi(a,b,u)>0$ for every $u\in U$. If $v=hu$ for some $h\in
G_{a,b}$ then $\phi(a,b,u)=\phi(a,b,v)$.
Let $g\in G$ be such that
$g(a,b)=(b,a)$. Then $g^{-1}G_{a,b}g=G_{a,b}$ and if $G$ is not
directed then for every $G_{a,b}$-invariant measurable subset $U$
of $\Lambda-\{a,b\}$ with $\mu(U)>0$ we have $\mu(U\cap gU)>0$.
Now for every $c\in \mathbb{R}$ and every $\epsilon >0$ the set
$U_{c,\epsilon}=\{u\mid \phi(a,b,u)\in (c-\epsilon,c+\epsilon)\}$
is $G_{a,b}$-invariant. By assumption, if $\mu(U_{c,\epsilon})>0$
then we have $\mu(gU_{c,\epsilon}\cap U_{c,\epsilon})>0$ as well
and therefore there is some $u\in U_{c,\epsilon}$ with $ \vert
\phi(a,b,gu)-c\vert <\epsilon$. Since $\phi$ is alternating we
obtain that
$\phi(a,b,gu)=-\phi(b,a,gu)=-\phi(ga,gb,gu)=-\phi(a,b,u)\in
(-c-\epsilon,-c+\epsilon)$ and consequently $\vert c\vert
<2\epsilon$. On the other hand, for $\mu$-almost every $u\in
\Lambda$ and every $\epsilon >0$ we have
$\mu(U_{\phi(a,b,u),\epsilon})>0$. From this we deduce that $\phi$
vanishes almost everywhere and hence
$H_{cb}^2(G,\mathbb{R})=\{0\}$. \qed
\bigskip

In the case that $G$ is a simple Lie group of non-compact
type and rank 1 we have $H_{cb}(G,\mathbb{R})
\not=\{0\}$ if and only if $G=SU(n,1)$ for some $n\geq 1$.
That this condition is necessary is immediate from Lemma 3.2
and the following observation.

\bigskip

{\bf Lemma 3.3:} {\it
A simple Lie group $G$ of non-compact type and rank one
is directed only if $G=SU(n,1)$ for some $n\geq 1$.
}

{\it Proof:} Let $G$ be a simple Lie group of non-compact type
and rank one. Then $G$ is the isometry group of a symmetric
space $X$ of non-compact type and negative curvature with
\emph{ideal boundary} $\partial X$.
The action of $G$ on $\partial X$ preserves the measure class
of the Lebesgue measure $\lambda$ and $(\partial X,\lambda)$
is a strong boundary for $G$.
The group $G$ acts
transitively on the space $A$
of pairs of distinct points in $\partial X$.
If $G=SO(n,1)$
for some $n\geq 3$ then $G$ acts transitively on the
space of triples of pairwise distinct points in $\partial X$ and
hence by the above, $G$ is not directed.

Now let $G=Sp(n,1)$ for some $n\geq 1$. Then for every pair
$(a,b)\in A$ there is a unique totally geodesic embedded
quaternionic line $L\subset X$ of constant curvature $-4$ whose
boundary $\partial L=S^3\subset\partial X$ contains $a$ and $b$.
The stabilizer $G_{a,b}$ of $(a,b)$ is contained in the stabilizer
$G_L$ of $L$ in $G$ which is conjugate to the quotient of the
group $Sp(1,1)\times Sp(n-1)<Sp(n,1)$ by its center. Moreover,
$G_L$ acts transitively on the space of triples of pairwise
distinct points in $\partial L$. In particular, for every $u\in
\partial L$ and every $g\in G$ with $g(a,b)=(b,a)$ the
$G_{a,b}$-orbit of $u$ coincides with the $G_{a,b}$-orbit of $gu$.
More generally, this is true for \emph{every} point $u\in \partial
X-\{a,b\}$. Namely, let $P:\partial X\to L$ be the shortest
distance projection. The subgroup $\tilde G_{Pu}$ of $G_L$ which
stabilizes $Pu$ is the quotient of the group $Sp(1)\times Sp(n-1)$
by its center. The factor subgroup $Sp(n-1)$ acts transitively on
$P^{-1}(Pu)$ while the orbit of $Pu$ under the group $G_{a,b}$
consists of the set of all points in $L$ whose distance to the
geodesic connecting $b$ to $a$ coincides with the distance of
$Pu$. As above we conclude that for every $u\in \partial X$ and
every $g\in G$ with $g(a,b)=(b,a)$ the $G_{a,b}$-orbit of $u$
coincides with the $G_{a,b}$-orbit of $gu$. In other words, $G$ is
not directed. In the same way we obtain that the exceptional Lie
group $F^4_{20}$ is not directed as well. \qed

\bigskip

Note that for every $n\geq 1$ the group $SU(n,1)$ admits a
non-trivial continuous second bounded cohomology class $\alpha\in
H_{cb}^2(SU(n,1),\mathbb{R})$ induced by the \emph{K\"ahler form}
$\omega$ of the complex hyperbolic space $X=SU(n,1)/S(U(n)U(1))$,
and this class generates $H_{cb}^2(SU(n,1),\mathbb{R})$.
The anti-symmetric cross ratio $[\,]$ on $\partial X$ defining the
class $\alpha$ is a continuous function on $\partial X$ with
values in $[-\pi,\pi]$. We have $\vert [a,b,c,d]\vert =\pi$ if and
only if $(a,c,b,d)$ or $(c,a,d,b)$ is an ordered quadruple of
points in the boundary of a complex line in $X$ (for a discussion
of this fact, references and applications to rigidity, see the
recent preprint of Burger and Iozzi \cite{BI05}).

\section{Rigidity of cocycles}

In this section we consider as before
a proper hyperbolic geodesic metric space $(X,d)$.
We equip the isometry group ${\rm Iso}(X)$ with
the compact open topology; with respect to this
topology, ${\rm Iso}(X)$ is a locally compact
$\sigma$-compact topological group.
Let $S$ be a standard Borel space and let
$\mu$ be a Borel probability
measure on $S$. Assume that $G$ is a locally compact
$\sigma$-compact topological group which
admits a measure preserving ergodic action on
a standard probability space
$(S,\mu)$. This action then defines a natural
continuous unitary representation of $G$ into
the Hilbert space $L^2(S,\mu)$ of square integrable
functions on $S$.
Let $\alpha:G\times S\to {\rm Iso}(X)$ be a cocycle,
i.e. $\alpha:G\times S\to {\rm Iso}(X)$ is a Borel
map which satisfies $\alpha(gh,x)=\alpha(g,hx)\alpha(h,x)$
for all $g,h\in G$ and $\mu$-almost every $x\in S$.
The cocycle is \emph{cohomologous} to a cocycle
$\beta:G\times S\to {\rm Iso}(X)$ if there is a measurable
function $\psi:S\to {\rm Iso}(X)$ such that
$\psi(gx)\alpha(g,x)=\beta(g,x)\psi(x)$ for all
$g\in G$ and almost all $x\in S$.
We use Theorem 4.1 of \cite{H05b} to show
(compare \cite{MS04}, \cite{MMS04}).

\bigskip

{\bf Lemma 4.1:} {\it Let $G$ be a locally compact
$\sigma$-compact group and let
$(S,\mu)$ be an ergodic $G$-probability space. Let $X$ be a proper
hyperbolic geodesic metric space and let $\alpha:G\times S\to {\rm
Iso}(X)$ be a cocycle. Let $H<{\rm Iso}(X)$ be a closed subgroup
with limit set $\Lambda$
and assume that $\alpha(G\times S)\subset H$ but that $\alpha$
is not cohomologous to a cocycle into
a proper subgroup of $H$.
If $H_{cb}^2(G,L^2(S,\mu))$ is finite dimensional
then either $H$ is elementary or $H$ acts transitively
on the complement of the diagonal in $\Lambda\times\Lambda$.}

{\it Proof:} Let $G$ be a locally
compact $\sigma$-compact group
which admits an ergodic measure preserving action on a
standard Borel probability space $(S,\mu)$.
Let $\alpha:G\times
S\to {\rm Iso}(X)$ be a cocycle and assume that $\alpha(G\times S)$ is
contained in a closed subgroup $H$ of ${\rm Iso}(X)$ and that
$\alpha$ is not equivalent to a cocycle with values in a proper
closed subgroup of $H$. Assume moreover that
$H$ is not elementary and does not
act transitively on the complement of the diagonal
in $\Lambda\times \Lambda$. By Theorem 4.1 of \cite{H05b},
$H_{cb}^2(H,\mathbb{R})$ is infinite
dimensional. More precisely, if $T\subset \Lambda^3$
is the space of triples of pairwise distinct points in
$\Lambda$ then
there is an infinite dimensional
vector space of
\emph{continuous} $H$-invariant cocycles on $T$, i.e.
continuous $H$-invariant maps $T\to \mathbb{R}$ which
satisfy the cocycle equation (\ref{triple}).

Let $(B,\nu)$ be a strong boundary for $G$. By Lemma 2.2 of
\cite{H05b} there is an $\alpha$-equivariant measurable map
$\phi:B\times S\to \Lambda$. This means that for almost every
$(w,\sigma)\in B\times S$ and every $g\in G$ we have
$\phi(gw,g\sigma)=\alpha(g,\sigma)\phi(w,\sigma)$. We claim that
the image of $\phi$ is dense in $\Lambda$. Namely, otherwise there
is a proper closed subset $A$ of $\Lambda$ which contains the
image of $\phi$. Then for every $g\in G$ and almost every $u\in S$
the element $\alpha(g,u)\in H$ stabilizes $A$, and $\alpha$ is
cohomologous to a cocycle with values in the intersection of $H$ with
the stabilizer of $A$.
Since $A$ is a closed proper subset of
$\Lambda$, the stabilizer of $A$ intersects $H$ in
a closed subgroup
$H$ of infinite index which contradicts our assumption on $\alpha$.

We now follow Monod and Shalom \cite{MS04}. Namely, choose a
continuous $H$-invariant bounded cocycle $\omega:T\to \mathbb{R}$.
Define a measurable map $\eta(\omega)$ from the space $\tilde B$
of triples of pairwise distinct points in $B$ into
$L^\infty(S,\mu)$ by
$\eta(x_1,x_2,x_3)(z)=\omega(\phi(x_1,z),\phi(x_2,z),\phi(x_3,z))$;
then $\eta$ can be viewed as a measurable $G$-equivariant bounded
cocycle with values in the Hilbert space of square integrable
function on $S$. Since $\omega$ is continuous and $\phi(B\times
S)$ is dense in $\Lambda$, this cocycle does not vanish
identically. As a consequence, the space of $G$-equivariant
measurable bounded cocycles $\tilde B\to L^2(S,\mu)$ is infinite
dimensional. On the other hand, $B$ is a strong boundary for $G$
and therefore $H_{b}^2(G,L^2(S,\mu))$ is isomorphic to the space
of measurable $G$-equivariant bounded cocycles $\tilde B\to
L^2(S,\mu)$ \cite{M}. In other words, $H_{cb}^2(G,L^2(S,\mu))$ is
infinite dimensional. This shows the lemma. \qed

\bigskip

Now consider a semi-simple Lie group $G$ with finite center and no
compact factors. Let $\Gamma<G$ be an irreducible lattice in $G$,
i.e. $\Gamma$ is a discrete subgroup of $G$ such that the volume
of $G/\Gamma$ is finite and that moroever if $G=G_1\times G_2$ is
a non-trivial product then the projection of $\Gamma$ to each
factor group $G_i$ $(i=1,2)$ is dense. The group $G$ acts on
$G/\Gamma$ preserving the projection of the Lebesgue measure
$\lambda$ on $G$.

Let $(S,\mu)$ be a standard probability space
with a measure preserving action of $\Gamma$. Then 
$\Gamma$ admits a measure preserving action on the product 
space $G\times S$. Since $\Gamma<G$ is a lattice, the
quotient space $(G\times S)/\Gamma$ can be
viewed as a bundle over $G/\Gamma$ with fibre $S$.
If $\Omega\subset G$ is a 
finite measure Borel fundamental domain
for the action of $\Gamma$ on $G$ then $\Omega\times S\subset
G\times S$ is a finite measure Borel fundamental domain for
the action of $\Gamma$ on $G\times S$ and 
up to normalization, the
product measure $\lambda\times \mu$ projects to a finite
measure $\nu$ on $(G\times S)/\Gamma$. The action of $G$
on $G\times S$ by left translation commutes with the
action of $\Gamma$ and hence it projects to an action on
$(G\times S)/\Gamma$ preserving the measure $\nu$ 
(see p.75 of \cite{Z}).   
We have.

\bigskip

{\bf Lemma 4.2:} {\it Let $\Gamma$ be an irreducible lattice in a
product $G=G_1\times G_2$ of two semi-simple non-compact Lie
groups $G_1,G_2$. If the action of $\Gamma$ on $(S,\mu)$ is mixing
then the induced action of $G_1$ on $((G\times S)/\Gamma,\nu)$
is ergodic.}

{\it Proof:} Let $\Gamma$ be an irreducible lattice in a product
$G=G_1\times G_2$ of two semi-simple non-compact Lie groups
$G_1,G_2$ with finite center. Let $(S,\mu)$ be a mixing
$\Gamma$-probability space. The induced action of $G$ on
$(G\times S)/\Gamma$ preserves the measure $\nu$
and restricts to an action of $G_1$. Since the measure $\nu$ is
finite, this action is ergodic if and only if every
$G_1$-invariant function $f\in L^2((G\times S)/\Gamma,\nu)$ is
constant.

We argue by contradiction and we assume that there is a
$G_1$-invariant function $f\in L^2((G\times S)/\Gamma,\nu)$ with
$\int \Vert f\Vert^2d\nu=1$ and zero mean $\int fd\nu=0$. 
Let $\tilde f$ be the lift of $f$ to a $\Gamma$-invariant function
on $G\times S$; then $\tilde f$ is locally square integrable
with respect to the product measure $\lambda\times \mu$.
For $x\in G$ define a function $\phi(x)$ on $S$ by
$\phi(x)(y)=\tilde f(x,y)$;
by Fubini's theorem we have $\phi(x)\in
L^2(S,\mu)$ for $\lambda$-almost every $x$, moreover the function
$x\to \int_S \phi(x) d\mu$ is measurable with respect to
$\lambda$. We claim that $\int_S\phi(x)d\mu=0$ and $\int_S\vert
\phi(x)\vert^2d\mu=1$ for $\lambda$-almost every $x\in G$.

For this recall that by definition of the action of $G_1$, the
functions $x\to \int_S\phi(x)d\mu$ and $x\to \int_S\vert
\phi(x)\vert^2d\mu$ are invariant under the action of $G_1$ on
$G/\Gamma$ since the action of $\Gamma$ on $(S,\mu)$ is measure
preserving. By Moore's ergodicity theorem the action of $G_1$ on
$G/\Gamma$ is ergodic \cite{Z} and hence these functions are
constant almost everywhere. If $\Omega\subset G$
is a Borel-fundamental domain for the action of $\Gamma$
on $G$ and if the Lebesgue measure $\lambda$ on $G$
is normalized in such a way that $\lambda(\Omega)=1$ then
$\int_{\Omega} \int_S\phi(x)d\mu
d\lambda=0$, $\int_{\Omega}\int_S \vert \phi(x)\vert^2d\mu d\lambda=1$ and
hence these constants equal $0,1$. In particular, if we denote
by ${\cal V}$ the orthogonal complement in $L^2(S,\mu)$ of the
constant functions then the assignment $x\to \phi(x)$ is a ${\cal
V}$-valued measurable function on $G/\Gamma$.

Let $L^{[1]}(G,{\cal V})^\Gamma$ be the space of all
$\Gamma$-invariant measurable maps $h:G\to {\cal V}$ for
which the function $\Vert h\Vert^2$ is locally integrable with respect to
the Lebesgue measure $\lambda$ on $G$. Note
that $L^{[1]}(G,{\cal V})^\Gamma$ naturally has the structure of a
separable Banach space. The group $G_1$ acts isometrically on 
$L^{[1]}(G,{\cal V})^\Gamma$ by left translation.
By our assumption, there
is a non-zero $G_1$-invariant vector $\phi\in 
L^{[1]}(G,{\cal V})^\Gamma$.
Then $\phi:G\to {\cal V}$ is measurable and essentially
bounded. Thus the restriction of $\phi$ 
to suitable compact sets of large measure
is continuous; in particular, we can find
a compact set $K=K_1\times K_2\subset G_1\times G_2$
of positive measure, Borel sets $D,E\subset S$ of positive measure
and a number $\epsilon >0$ such that the following is satisfied.
\begin{enumerate}
\item
$\phi(z)(u)\geq \epsilon$ for every $z\in K, u\in D$.
\item $\phi(z)(v)
\leq 0$ for every $z\in K,v\in E$.
\end{enumerate}

Since $\phi$ is invariant
under the left action of $G_1$, for any $z_0\in K_1,z_2\in K_2$ 
and $z_1\in G_1$, 
$u\in S$ we have
$\phi(z_1,z_2)(u)=\phi(z_0,z_2)(u)$. 
On the other hand, 
$\phi$ is also invariant under the right action of $\Gamma$ and hence
$\phi(z_1,z_2)(u)=\phi((z_1,z_2)\eta)(\eta u)$
for every $\eta\in \Gamma$. In particular, if $(z_1,z_2)\eta\in K$
then we deduce that $\mu(\eta D\cap E)=0$.

The right action of $\Gamma$ on $G_1\backslash G$ is ergodic and
therefore there are infinitely many elements $\eta\in \Gamma$ such
that the measure of $(G_1\times K_2)\eta\cap K$ is positive. 
Moreover,
the action of $\Gamma$ on
$(S,\mu)$ is mixing by assumption and hence for every $\eta\in
\Gamma$ which is sufficiently far away from the identity the image
of $D$ under $\eta$ intersects $E$ in a set of positive measure.
In particular, there is some $\eta\in \Gamma$ with
$\lambda((G_1\times K_2)\eta\cap K)>0$ and $\mu(\eta D\cap E)>0$
which contradicts the above consideration.
Thus the vector space of $G_1$-invariant maps in
$L^{[1]}(G,{\cal V})^\Gamma$ is trivial and the action of $G_1$
on $(G\times S)/\Gamma$ is ergodic.
\qed

\bigskip

The following corollary completes the proof of
Theorem A from the introduction and follows as in
\cite{MS04} from Lemma 4.2, the
rigidity results for bounded cohomology
of Burger and Monod \cite{BM99,BM02} and
the work of Zimmer \cite{Z}.

\bigskip

{\bf Corollary 4.3:} {\it Let $G$ be a semi-simple Lie group with
finite center, no compact factors and of rank at least $2$. Let
$\Gamma<G$ be an irreducible lattice and let $(S,\mu)$ be a mixing
$\Gamma$-space. Let $X$ be a proper hyperbolic geodesic metric
space and let $\alpha:\Gamma\times S\to {\rm Iso}(X)$ be a
cocycle; then either $\alpha$ is cohomologous to a cocycle into an
elementary subgroup of ${\rm Iso}(X)$ or there is a closed
subgroup $H$ of ${\rm Iso}(X)$ which is a compact extension of a
simple Lie group $L$ of rank one and there is a surjective
homomorphism $G\to L$.}

{\it Proof:} Let $G$ be a semi-simple Lie group of non-compact
type with finite center and of rank at least 2 and let $\Gamma$ be
an irreducible lattice in $G$. Let $(S,\mu)$ be a mixing
$\Gamma$-space with invariant Borel probability measure $\mu$. Let
$\alpha:\Gamma\times S\to {\rm Iso}(X)$ be a cocycle into the
isometry group of a proper hyperbolic geodesic metric space $X$.
Let $H<{\rm Iso}(X)$ be a closed subgroup such that $\alpha$ is
cohomologous to an $H$-valued cocycle but not to a cocycle with
values in a proper subgroup of $H$. We may assume that
$\alpha(\Gamma\times S) \subset H$. Assume that $H$ is not
elementary, with limit set $\Lambda$. By Theorem A in \cite{H05b},
there is a continuous $L^2(H)$-valued bounded cocycle
$\phi:\Lambda^3\to L^2(H)$.

Let $\Omega\subset G$ be a Borel fundamental domain for
the action of $\Gamma$ on $G$. Let $\nu$ be the $G$-invariant
Borel probability measure on $(G\times S)/\Gamma\sim \Omega\times S$. 
We obtain
a $\nu$-measurable function $\beta:G\times (G\times S)/\Gamma\to 
H$ as follows. For $z\in \Omega$ and $g\in G$
let $\eta(g,z)\in \Gamma$ be the unique
element such that $gz\in \eta(g,z)\Omega$ and define
$\beta(g,(z,\sigma))=\alpha(\eta(g,z),\sigma)$. By construction,
$\eta$ satisfies the cocycle equation for the action of $G$ on
$(G\times S)/\Gamma$. Let
$(B,\mu)$ be a strong boundary for $G$; we may assume that $B$ is
also a strong boundary for $\Gamma$. By Lemma 4.2 the action
of $G$ on $(G\times S)/\Gamma$ is ergodic and hence Lemma 2.2 of
\cite{H05b} show that we can find a measurable $\beta$-equivariant
F\"urstenberg map $\psi_0:(G\times S)/\Gamma\times B\to \Lambda$ and
hence a $\beta$-equivariant map $\psi:(G\times S)/\Gamma\times
B^3\to \Lambda^3$ defined by $\psi(x,a,b,c)=
(\psi_0(x,a),\psi_0(x,b), \psi_0(x,c))$ $(x\in (G\times S)/\Gamma)$.
Let $L^{[2]}((G\times S)/\Gamma,L^2(H))$ be the space of all
measurable maps $(G\times S)/\Gamma\to L^2(H)$ with the additional
property that for each such map $\phi$ the function $x\to
\Vert\phi\Vert$ is square integrable on $(G\times S)/\Gamma$. Then
$L^{[2]}((G\times S)/\Gamma,L^2(H))$ has a natural structure of a
separable Hilbert space, and the group $G$ acts on
$L^{[2]}((G\times S)/\Gamma,L^2(H))$ as a group of isometries. In
other words, $L^{[2]}((G\times S)/\Gamma,L^2(H))$ is a Hilbert
module for $G$ and the cocycle $\phi$ can be composed with the map
$\psi$ to a $\beta$-invariant measurable bounded map $B^3\to
L^{[2]}((G\times S)/\Gamma,L^2(H))^3$. Since $B$ is a strong
boundary for $G$ we conclude that this map defines a nontrivial
cohomology class in $H_{cb}^2(G,L^{[2]}((G\times S)/\Gamma,
L^2(H)))$.

Now if $G$ is simple then the results of Monod
and Shalom \cite{MS04} show that there is a $\beta$-equivariant map
$(G\times S)/\Gamma\to L^2(H)$. Since the action of $G$ on
$(G\times S)/\Gamma$ is ergodic, by the cocycle reduction lemma of
Zimmer \cite{Z} the cocycle $\beta$ and hence $\alpha$ is
cohomologous to a cocycle into a compact subgroup of $H$ which is
a contradiction. On the other hand, if $G=G_1\times G_2$ for
semi-simple Lie groups $G_1,G_2$ with finite center and without
compact factors, then the results of Burger and Monod
\cite{BM99,BM02} show that via possibly exchanging $G_1$ and $G_2$
we may assume that there is an equivariant map $(G\times S)/\Gamma
\to L^2(H)$ for the restriction of $\beta$ to $G_1\times
(G\times S)/\Gamma$, viewed as a cocycle for $G_1$. By Lemma 4.2 the
action of $G_1$ on $(G\times S)/\Gamma$ is ergodic and therefore the
cocycle reduction lemma of Zimmer \cite{Z} shows that the
restriction of $\beta$ to $G_1$ is equivalent to a cocycle into a
compact subgroup of $H$. We now follow the proof of Theorem 1.2 of
\cite{MS04} and find a minimal such compact subgroup $K$ of $H$.
The cocycle $\beta$ and hence $\alpha$ is cohomologous to a cocycle
into the normalizer of $K$ in $H$ which then coincides with $H$ by
our assumption of $H$. Moreover, there is a continuous
homomorphism of $G$ onto $H/K$. Since $G$ is connected, the image
of $G$ under this homomorphism is connected as well and hence by
Theorem A of \cite{H05b}, $H/K$ is a simple Lie group of rank one.
This completes the proof of the corollary. \qed

\bigskip

\noindent
MATHEMATISCHES INSTITUT DER UNIVERSIT\"AT BONN\\
BERINGSTRA\SS{}E 1\\
D-53115 BONN\\

\smallskip

\noindent
e-mail: ursula@math.uni-bonn.de

\end{document}